**ORIGINAL ARTICLE**

# Convergence results for an averaged LQR problem with applications to reinforcement learning

**Andrea Pesare[1]** · **Michele Palladino[2]** · **Maurizio Falcone[1]**

## Abstract

In this paper, we will deal with a linear quadratic optimal control problem with unknown dynamics. As a modeling assumption, we will suppose that the knowledge that an agent has on the current system is represented by a probability distribution $\pi$ on the space of matrices. Furthermore, we will assume that such a probability measure is opportunely updated to take into account the increased experience that the agent obtains while exploring the environment, approximating with increasing accuracy the underlying dynamics. Under these assumptions, we will show that the optimal control obtained by solving the "average" linear quadratic optimal control problem with respect to a certain $\pi$ converges to the optimal control driven related to the linear quadratic optimal control problem governed by the actual, underlying dynamics. This approach is closely related to model-based reinforcement learning algorithms where prior and posterior probability distributions describing the knowledge on the uncertain system are recursively updated. In the last section, we will show a numerical test that confirms the theoretical results.

✉ Maurizio Falcone
falcone@mat.uniroma1.it

Andrea Pesare
pesare@mat.uniroma1.it

Michele Palladino
michele.palladino@gssi.it

1 Dipartimento di Matematica, Sapienza - Università di Roma, Rome, Italy

2 Gran Sasso Science Institute - GSSI, L'Aquila, Italy

## 1 Introduction

Reinforcement learning (RL) is one of the three basic machine learning paradigms, together with supervised learning and unsupervised learning. In RL, an agent interacts with a partially unknown environment, aiming at finding a policy, which optimizes the measure of a certain long-term performance [25]. The connection between RL and optimal control theory was already identified in the past [26]. Nowadays, several research questions and tools coming from the RL literature are influencing the optimal control field and vice versa [23].

A natural setting in RL consists in considering a Markov decision process over state/action pairs varying on a discrete-time set. The discrete-time problem setting provides an excellent framework to develop methods and algorithms, which, however, often underlies a continuous-time structure. For this reason, in particular in the control system engineering field, significant attention has been recently given to continuous-time RL [8,14,15,17].

Both in discrete- and continuous-time problem settings, one can consider two main RL philosophies: The first one, called *model-based*, usually concerns the reconstruction of a model from the data trying to mimic the unknown environment. That model is then used to plan and to compute a suboptimal policy. The second RL philosophy, called *model-free*, employs a direct approximation of the value function and/or a policy based on a dynamic-programming-like algorithm, without using a model to simulate the unknown environment. An excellent overview of the two approaches can be found in [25].

About ten years ago, PILCO was introduced [6,7], an innovative and disruptive method, from which many subsequent model-based RL methods have been inspired. Rather than exploiting the data to construct a dynamics approximating the partially known environment, PILCO makes use of them to construct a probability distribution (more precisely, a Gaussian process) on a class of dynamical systems. At each iteration, this distribution is updated to fit the data set. After the model update, PILCO takes the *policy improvement* step which boils down to solving an averaged optimal control problem, where the averaging distribution is the one extrapolated by the data at the previous experiments. That approach has the advantage of considerably reducing the model bias, one of the main shortcomings of model-based RL [1]. A general, rigorous framework capturing PILCO as well as other Bayesian model-based RL approaches (see, e.g., [3,4,10–12,29]) has been developed in [18,19]. In particular, it is important to mention that the framework developed in [18] is closely related to the averaging control framework and Riemann–Stieltjes optimal control [2,16,21,24,30].

The aim of this paper is to provide a stricter link between PILCO [7] and the framework introduced in [18]. We will concentrate on a specific physical system, driven by a linear quadratic regulator (LQR) optimal control problem, namely

$$\begin{cases} \text{minimize } \left\{ \frac{1}{2} \int_0^T \left( x(t)^T Q x(t) + u(t)^T R u(t) \right) dt + \frac{1}{2} x(T)^T Q_f x(T) \right\} \\ \text{over } u : [0, T] \to \mathbb{R}^m \text{ measurable such that} \\ \dot{x}(t) = \hat{A} x(t) + B u(t), \qquad t \in [0, T] \\ x(0) = x_0 \end{cases} \tag{1}$$

where $B$, $Q$, $R$, $Q_f$ are given, known matrices and $x_0 \in \mathbb{R}^n$ is a given vector. However, in this context we will assume that the physical system $\hat{A}$ is *unknown* and the knowledge we have about $\hat{A}$ is merely represented by a probability distribution $\pi$ constructed over a set of matrices $\mathcal{A}$ (with $\hat{A} \in \mathcal{A}$) by using the data available from the physical system. This situation is in accordance with the PILCO setting, which is "not focusing on a single dynamics model", but makes use of "a probabilistic dynamics model, a distribution over all plausible dynamics models that could have generated the observed experience" ( [6], pg. 34). Such a modeling setting allows us to define the averaged optimal control problem

$$
\begin{cases}
\text{minimize } \left\{ \int_{\mathcal{A}} \int_0^T \frac{1}{2} \left( x_A(t)^T Q x_A(t) + u(t)^T R u(t) \right) dt \right. \\
\qquad\qquad\qquad\qquad \left. + \frac{1}{2} x_A(T)^T Q_f x_A(T) d\pi(A) \right\} \\
\text{over } u : [0, T] \to \mathbb{R}^m \text{ measurable such that} \\
\dot{x}_A(t) = A x_A(t) + B u(t), \quad A \in \mathcal{A}, \quad t \in [0, T], \\
x_A(0) = x_0
\end{cases}
\tag{2}
$$

If indeed the real physical system is driven by the equation

$$
\dot{x}(t) = \hat{A} x(t) + B u(t), \qquad t \in [0, T],
$$

for a certain matrix $\hat{A}$, then it is reasonable to expect that an increase of the experience will produce a more accurate distribution $\pi$ over $\mathcal{A}$. This fact can be translated into the assumption that the probability distribution is "close" (in a precise sense that will be specified in the sequel) to $\delta_{\hat{A}}$, when enough experience of the environment (here represented by $\hat{A}$) is gained.

We would like to stress that our goal is not to propose a new algorithm to find an optimal policy but to consider a class of existing algorithms and motivate their good performances. In particular, the paper aims to provide an insight into the convergence of Bayesian-like RL algorithms in which a recursive construction of probability measures is carried out. (Further considerations on the connection with RL are given in Remarks 3.4 and 5.5.) Here, by "convergence", we mean convergence of the optimal policy obtained by estimating the underlying dynamics using data from the real system (the one constructed in the so-called *policy improvement* step) toward the optimal policy obtained by solving problem (1). More precisely, the questions we will tackle in this paper are:

(1) Is the value function related to the optimal control (1) close to the value function associated with (2) when $\pi$ is close to $\delta_{\hat{A}}$ w.r.t. the Wasserstein distance (see (9) for the formal definition)?
(2) Under the same assumptions over $\pi$, is the optimal control of (2) close to the optimal control of (1)?

For both question (1) and question (2), we will provide positive answers. It is worth noticing that, in control theory, it is very uncommon to have a positive answer to question (2), even when one has a positive answer to question (1).

The paper is organized as follows. Section 2 introduces the basic notations that we will use throughout the paper. In Sect. 3, we state the problem formulation and we study the basic properties of the LQ optimal control problem (2). Then, in Sect. 4, we also derive a Pontryagin's maximum principle for problem (2), refining some results in [2]. In Sect. 5, we state and prove the main results of the paper, providing positive answers to question (1) and question (2). In Sect. 6, we strengthen the results presented in Sect. 5 in the case in which one is dealing with a discrete probability measure $\pi$. That result is further stressed in Sect. 7, where we present and analyze a numerical example. Finally, we present the conclusion and discuss some future directions and open questions.

## 2 Preliminaries and notations

In this section, we will recall some useful notations and concepts which will be used throughout the paper. For vectors $v \in \mathbb{R}^n$, $|v|$ denotes the Euclidean norm; we use $\mathbb{B}_n(x, r)$ to denote the open ball in $\mathbb{R}^n$ centered at $x \in \mathbb{R}^n$ and of radius $r > 0$. We will use $\mathcal{L}$ to denote the Lebesgue $\sigma$-algebra and $\mathcal{B}_X$ to denote the collection of all Borel sets on a given topological space $X$.

In the following, $T$ will be our fixed time horizon. Given a time $s \in [0, T]$, we denote by $C([s, T]; \mathbb{R}^n)$ the space of continuous functions $x \colon [s, T] \to \mathbb{R}^n$. For functions $x(\cdot) \in C([s, T]; \mathbb{R}^n)$, $\|x(\cdot)\|_\infty$ or $\|x\|_\infty$ denotes the sup norm. For continuous functions $y(\cdot) \in C(\mathbb{R}^n; \mathbb{R})$ and for a compact set $K \subset \mathbb{R}^n$, we also define the sup norm restricted on the compact set $K$ as $\|y(\cdot)\|_{\infty, K} := \sup_{x \in K} |y(x)|$.

Moreover, given $p \in [1, \infty)$, we define the spaces of a.e. defined functions

$$L^p([s, T]; \mathbb{R}^n) := \left\{ g \colon [s, T] \to \mathbb{R}^n \,\middle|\, g \text{ meas. and } \int_s^T |g(t)|^p \, dt < \infty \right\} \quad (3)$$

and

$$W^{1,p}([s, T]; \mathbb{R}^n) := \left\{ g \in L^p([s, T]; \mathbb{R}^n) \,\middle|\, \frac{dg}{dt} \in L^p([s, T]; \mathbb{R}^n) \right\}. \quad (4)$$

For $g \in L^p([s, T]; \mathbb{R}^n)$, the $L^p$-norm $\|g(\cdot)\|_p$ or $\|g\|_p$ is defined by

$$\|g\|_p := \left( \int_s^T |g(t)|^p \, dt \right)^{\frac{1}{p}} \quad (5)$$

and for $g \in W^{1,p}([s, T]; \mathbb{R}^n)$, we define

$$\|g\|_{W^{1,p}} := \|g\|_p + \left\|\frac{dg}{dt}\right\|_p. \quad (6)$$

Let us denote by $\mathbb{M}_{m \times n}$ the space of real matrices with $m$ rows and $n$ columns. For square matrices $A \in \mathbb{M}_{n \times n}$, we consider the $2$−norm

$$\|A\|_2 := \sup \left\{ x^T A y \colon \ x, y \in \mathbb{R}^n, \ |x| = |y| = 1 \right\}. \quad (7)$$

Given two matrices $A, A' \in \mathbb{M}_{n \times n}$, $d_2(A, A') := \left\| A - A' \right\|_2$ denotes the distance between the two matrices induced by the $2-$norm.

For a generic metric space $(X, d)$, $\mathcal{M}(X)$ will denote the space of measures on $X$, equipped with the weak-$\star$ topology, according to which $\pi^N \overset{*}{\rightharpoonup} \pi$ if and only if

$$\int_X f(A) d\pi^N(A) \to \int_X f(A) d\pi(A), \qquad \forall f \in C_b(X), \tag{8}$$

where $C_b(X)$ denote all bounded continuous functions $f : X \to \mathbb{R}$. When the space $X$ is compact, the weak-$\star$ topology is metrized by the Wasserstein distance (see, e.g., [27])

$$W_1(\pi, \pi') := \inf_{\gamma \in \Gamma(\pi, \pi')} \int_{X \times X} d(x, y) \, d\gamma(x, y) , \tag{9}$$

where $\Gamma(\pi, \pi')$ is the collection of all probability measures on $X \times X$ with marginals $\pi$ and $\pi'$ on the first and second factors, respectively.

Given a probability space $(\Omega, \mathcal{F}, \pi)$ and a random variable $Y$ on $\Omega$, $\mathbb{E}_\pi[Y]$ denotes the expected value of $Y$ with respect to $\pi$.

## 3 Problem statements and preliminary results

We begin our discussion considering two LQ optimal control problems. We will see in the sequel how the two problems are connected:

*Problem A: the LQR problem.*

Let us consider the classical LQR problem with finite horizon, which we will refer to as *Problem A*:

$$\begin{cases} \text{minimize } J_s[u] \\ \text{over } (x, u)(\cdot) \text{ such that } u \in \mathcal{U} \text{ and} \\ \dot{x}(t) = \hat{A}x(t) + Bu(t), \quad t \in [s, T], \\ x(s) = x_0 \end{cases} \tag{10}$$

where $s \in [0, T]$, $x_0 \in \mathbb{R}^n$, $\mathcal{U} := \{u : [s, T] \to \mathbb{R}^m \text{ Lebesgue measurable}\}$ and

$$J_s[u] := \frac{1}{2} \int_s^T \left( x(t)^T Q x(t) + u(t)^T R u(t) \right) dt + \frac{1}{2} x(T)^T Q_f x(T) . \tag{11}$$

The pair $(x, u)(\cdot)$ such that $u \in \mathcal{U}$ and $x(\cdot)$ is the solution of the Cauchy problem

$$\begin{cases} \dot{x}(t) = \hat{A}x(t) + Bu(t) \qquad t \in [s, T] \\ x(s) = x_0 . \end{cases} \tag{12}$$

is called *admissible process for Problem A.*

Let us define the *value function* $V : [0, T] \times \mathbb{R}^n \to \mathbb{R}$ for Problem A as

$$V(s, x_0) := \inf_{u \in \mathcal{U}} J_s[u]. \tag{13}$$

We shall say that $(\bar{x}, \bar{u})(\cdot)$ is an *optimal process for Problem A* if

$$J_s[\bar{u}] \leq J_s[u] \tag{14}$$

for any other admissible process $(x, u)(\cdot)$ of Problem A. In this case, $\bar{u}$ will be denoted as *optimal control* for Problem A.

*Problem B: an LQR problem with unknown dynamics.* Let us now introduce an optimal control that does not require the exact knowledge of the matrix $\hat{A}$, but merely a probability distribution defined on a compact space of matrices $\mathcal{A}$ containing $\hat{A}$. For each $s \in [0, T]$, $x_0 \in \mathbb{R}^n$ and $\pi \in \mathcal{M}(\mathcal{A})$, consider the following optimal control problem, which we will refer to as *Problem B*:

$$\begin{cases} \text{minimize } J_{s,\pi}[u] \\ \text{over } \{(x_A, u)(\cdot) : \ A \in \mathcal{A}\} \text{ such that } u \in \mathcal{U} \text{ and} \\ \dot{x}_A(t) = A x_A(t) + B u(t), \quad A \in \mathcal{A}, \quad t \in [s, T], \\ x_A(s) = x_0 \quad A \in \mathcal{A}, \end{cases} \tag{15}$$

where $\mathcal{U} := \{u : [s, T] \to \mathbb{R}^m \text{ Lebesgue measurable}\}$ and

$$\begin{aligned} J_{s,\pi}[u] &:= \mathbb{E}_\pi \left[ \frac{1}{2} \int_s^T \left( x_A(t)^T Q x_A(t) + u(t)^T R u(t) \right) dt + \frac{1}{2} x_A(T)^T Q_f x_A(T) \right] \\ &= \int_{\mathcal{A}} \left[ \frac{1}{2} \int_s^T \left( x_A(t)^T Q x_A(t) + u(t)^T R u(t) \right) dt + \frac{1}{2} x_A(T)^T Q_f x_A(T) \right] d\pi(A). \end{aligned} \tag{16}$$

**Remark 3.1** Sometimes we will denote this problem as Problem $B_\pi$, to stress its dependency on the probability distribution $\pi$.

The definition of the value function $V_\pi : [0, T] \times \mathbb{R}^n \to \mathbb{R}$ for Problem B is

$$V_\pi(s, x_0) := \inf_{u \in \mathcal{U}} J_{s,\pi}[u]. \tag{17}$$

The collection $\{(x_A, u)(\cdot) : \ A \in \mathcal{A}\}$ such that $u \in \mathcal{U}$ and, for each $A \in \mathcal{A}$, $x_A(\cdot)$ is the solution of the Cauchy problem

$$\begin{cases} \dot{x}_A(t) = A x_A(t) + B u(t) \qquad t \in [s, T] \\ x_A(s) = x_0 \end{cases} \tag{18}$$

is called *admissible process for Problem B*. Note that the initial condition is the same for every $A$. The admissible process $\{(\bar{x}_A, \bar{u})(\cdot) : A \in \mathcal{A}\}$ is *optimal for Problem B* if

$$J_{s,\pi}[\bar{u}] \leq J_{s,\pi}[u] \tag{19}$$

for any other admissible process $\{(x_A, u)(\cdot) : A \in \mathcal{A}\}$ of Problem B. In this case, $\bar{u}$ will be denoted as *optimal control for Problem B*.

Throughout the whole paper, the following *Standing Hypothesis* will be imposed both for Problem A and for Problem B:

**(SH):** $Q, Q_f \in \mathbb{M}_{n \times n}$ are symmetric and semipositive definite, and $R \in \mathbb{M}_{m \times m}$ is symmetric and positive definite.

**Remark 3.2** As we will show in Lemma 3.6, Problem B admits a unique optimal process when the assumption **(SH)** holds true. Furthermore, it is interesting to observe that Problem A can be regarded as a particular case of Problem B when one chooses $\pi = \delta_{\hat{A}}$, namely when $\pi$ is a Dirac delta concentrated at $\hat{A}$:

$$\text{Problem B}_{\delta_{\hat{A}}} \hookrightarrow \text{Problem A}$$

**Remark 3.3** In our framework, we consider a probability distribution merely on the matrix $A$ of the dynamics and not on the matrix $B$. It is not hard to check that the arguments proposed and the results achieved in the next sections hold as well in an extended framework where the matrix $B$ is possibly unknown. However, we preferred to consider a simpler case to keep the overall presentation as clear as possible. Moreover, it seems reasonable to assume that the agent does not know how the environment works (matrix $A$), whereas being aware of how the control affects the system (matrix $B$).

**Remark 3.4** (*Nonlinear systems and connection with RL*) We defined $\pi$ as a probability distribution on a space of matrices $\mathcal{A}$. In a higher perspective, we can identify $\mathcal{A}$ with a class of *linear* dynamical systems, namely

$$\widetilde{\mathcal{A}} := \{(x, u) \mapsto Ax + Bu \,|\, A \in \mathcal{A}\},$$

and see $\pi$ as a probability distribution on $\widetilde{\mathcal{A}}$ as well.

More generally, one could consider a nonlinear dynamical system

$$\begin{cases} \dot{x}(t) = f(x(t), u(t)) & t \in [s, T] \\ x(s) = x_0, \end{cases} \tag{20}$$

assuming that the $f$ is unknown and only a probability distribution $\pi$ on a space $X$ of possible dynamics is available. In a similar way to how we did, given a cost functional

$$J_s[u] := \int_s^T \ell(x(t), u(t)) dt + h(x(T)), \tag{21}$$

one could define a Problem B for the nonlinear problem (see [18,22]).

A concrete example of *nonlinear* Problem B related to RL is presented in PILCO, by Deisenroth and Rasmussen [6,7]. PILCO aims at approximating the actual dynamics with a Gaussian process (GP); the parameters of the GP are tuned according to the data gathered by an agent while interacting with the environment. Providing a probability distribution of the dynamics permits to obtain an interval of confidence of the actual dynamics at each point, rather than a pointwise estimate of the actual dynamics. We stress that a GP is a probability distribution on the set $X$ of bounded, continuous functions, and thus, it constitutes a perfect candidate to play the role of $\pi$ in Problem B. In the *policy improvement* step, the optimal policy is then computed minimizing the averaged cost over all possible realizations of the GP, in a similar way to how we defined the cost functional in (16).

The framework presented in this paper incorporates the approach presented in PILCO, as well as other probabilistic model-based RL methods [29], when the dynamics is linear and the cost is quadratic. We will make further considerations about the link with PILCO in Remark 5.5.

### 3.1 Preliminary results for Problem B

Let us start showing a series of basic results on the existence and the regularity of trajectories and optimal controls for the system we are considering.

In this section, we assume that $\pi$ is a probability distribution on a compact set of matrices $\mathcal{A} \subset \mathbb{M}_{n \times n}$. Being $\mathcal{A}$ bounded, there exists a constant $C_A$ such that $\|A\|_2 \leq C_A, \ \forall A \in \mathcal{A}$.

For a given matrix $A \in \mathcal{A}$ and an admissible control $u \in \mathcal{U}$, the notation $x_A(t; u)$ will denote the solution of (18) relative to $A$ and $u$; sometimes, when it is not ambiguous, we could omit the dependency on the control $u$ and write only $x_A(t)$.

**Lemma 3.5** (Boundness and continuity of trajectories) *Let us consider the dynamical system in* (18). *The following results hold:*

*(i) For each $u \in L^1([s, T]; \mathbb{R}^m)$, the trajectory $x_A(\cdot; u)$ is uniformly bounded for all $A \in \mathcal{A}$:*

$$|x_A(t; u)| \leq C_x^u \quad \forall t \in [s, T], \ \forall A \in \mathcal{A},$$

*where $C_x^u$ is a constant which depends on $u$ and on $x_0$.*
*(ii) For each $u \in L^1([s, T]; \mathbb{R}^m)$, the map $A \mapsto x_A(\cdot; u)$ is continuous;*
*(iii) For each $A \in \mathcal{A}$, the map $u \mapsto x_A(\cdot; u)$ is continuous for $u \in L^1([s, T]; \mathbb{R}^m)$.*

**Proof** (i) Recall that $x_A(t)$ satisfies the relation

$$x_A(t) = x_0 + \int_s^t A x_A(\tau) + B u(\tau) \, d\tau, \quad \forall t \in [s, T].$$

for each $A \in \mathcal{A}$. Then,

$$|x_A(t)| \leq |x_0| + \int_s^t \|A\|_2 \, |x_A(\tau)| \, d\tau + \|B\|_2 \int_s^t |u(\tau)| \, d\tau \quad \forall \, t \in [s, T],$$

so by the Grönwall lemma (see, e.g., Lemma 2.4.4 on [28]) we get

$$
\begin{aligned}
|x_A(t)| &\leq \left( |x_0| + \|B\|_2 \int_s^t |u(\tau)| \, d\tau \right) e^{\|A\|_2(t-s)} \\
&\leq \left( |x_0| + \|B\|_2 \int_s^t |u(\tau)| \, d\tau \right) e^{C_A T} =: C_x^u
\end{aligned}
\tag{22}
$$

for every $A \in \mathcal{A}$ and $t \in [s, T]$. This shows condition $(i)$.

(ii) Fix a control $u \in L^1([s, T]; \mathbb{R}^m)$ and consider the trajectories solutions of (18) relative to two different matrices $A$ and $A'$. If we define

$$z(t) := x_A(t) - x_{A'}(t) \,,$$

then $z$ solves the following differential system:

$$
\begin{cases}
\dot{z}(t) = A x_A(t) - A' x_{A'}(t) & t \in [s, T] \\
z(s) = 0 \,.
\end{cases}
$$

Notice that we can rewrite the right-hand side as

$$A x_A(t) - A' x_{A'}(t) = A z(t) + (A - A') x_{A'}(t) \,,$$

and thus, we can give the estimate

$$|\dot{z}(t)| \leq C_A |z(t)| + \|A - A'\|_2 \, C_x^u \,.$$

Applying again the Grönwall lemma on $z$, we get

$$|z(t)| \leq \|A - A'\|_2 \, C_x^u \int_s^t e^{\int_\tau^t C_A d\sigma} d\tau \leq T C_x^u e^{C_A T} \, \|A - A'\|_2 \quad \forall \, t \in [s, T] \,,$$

which implies the continuity of the map $A \mapsto x_A(\cdot)$.

(iii) Fix a matrix $A \in \mathcal{A}$ and consider the trajectories relative to two different controls $u, u' \in \mathcal{U}$. In a similar way as in (ii), we define

$$z(t) := x_A(t; u) - x_A(t; u') \,,$$

which solves the ODE system

$$
\begin{cases}
\dot{z}(t) = A z(t) + B \big( u(t) - u'(t) \big) & t \in [s, T] \\
z(s) = 0 \,,
\end{cases}
$$

and, by Grönwall's lemma, we get

$$|z(t)| \leq \|B\|_2 \int_s^t e^{\int_\tau^t C_A d\sigma} |u(\tau) - u'(\tau)| d\tau \leq e^{C_A T} \|B\|_2 \|u - u'\|_1$$

for all $t \in [s, T]$. The last inequality gives the continuity with respect to $u \in L^1([s, T]; \mathbb{R}^m)$.

$\square$

The following result guarantees that the minimization problem (15) is well posed:

**Lemma 3.6** (Existence, uniqueness and upper bounds of the optimal control) *Let the assumption (SH) hold true. Given a $\pi \in \mathcal{M}(\mathcal{A})$, Problem B (15) admits a unique minimizer $\{(\bar{x}_A, \bar{u})(\cdot) : A \in \mathcal{A}\}$, satisfying the following upper bound:*

$$\left( \int_s^T |\bar{u}(t)|^2 \, dt \right)^{1/2} \leq \overline{C}_u \,, \tag{23}$$

*where $\overline{C}_u$ does not depend on $\pi$ and is defined as*

$$\overline{C}_u := \sqrt{\frac{1}{r_1} \left( T \|Q\|_2 + \|Q_f\|_2 \right) |x_0|^2 e^{2 C_A T}} \,,$$

*where $r_1$ is the smallest eigenvalue of the matrix R.*

**Proof** Consider a minimizing sequence $u^k \in \mathcal{U}$ of the cost functional $J_{s,\pi}$ defined in (16), satisfying $J_{s,\pi}[u^k] \to \inf_{u \in \mathcal{U}} J_{s,\pi}[u]$, $\inf_{u \in \mathcal{U}} J_{s,\pi}[u] \leq J_{s,\pi}[u^k]$, and the related minimizing process $\{(x_A^k, u^k)(\cdot) : A \in \mathcal{A}\}$. Set $\varepsilon_k := J_{s,\pi}[u^k] - \inf_{u \in \mathcal{U}} J_{s,\pi}[u] \geq 0$. It is not restrictive to assume that $\varepsilon_k < 1$ for all $k \in \mathbb{N}$.

Let us consider the system (18) when the control is $u^0 \equiv 0$. The process $\{(x_A^0, u^0)(\cdot) : A \in \mathcal{A}\})$ is solution of

$$\begin{cases} \dot{x}_A^0(t) = A x_A^0(t) & t \in [s, T] \\ x_A^0(s) = x_0 \end{cases}$$

for every $A \in \mathcal{A}$. Hence,

$$|x_A^0(t)| = |e^{At} x_0| \leq |x_0| e^{\|A\|_2 T} \quad \forall t \in [0, T], \quad \forall A \in \mathcal{A} \,.$$

Clearly, the cost achieved by the control $u^0$ can be estimated as follows:

$$J_{s,\pi}[u^0] = \mathbb{E}_\pi \left[ \frac{1}{2} \int_s^T x_A^0(t)^T Q x_A^0(t) \, dt + \frac{1}{2} x_A^0(T)^T Q_f x_A^0(T) \right]$$

$$\leq \mathbb{E}_\pi \left[ \frac{1}{2} T \|Q\|_2 |x_0|^2 e^{2\|A\|_2 T} + \frac{1}{2} \|Q_f\|_2 |x_0|^2 e^{2\|A\|_2 T} \right]$$

$$\leq \frac{1}{2} \left( T \left\| Q \right\|_2 + \left\| Q_f \right\|_2 \right) |x_0|^2 e^{2C_A T} .$$

Furthermore, it follows from the construction of the minimizing sequence that

$$J_{s,\pi}[u^k] \leq J_{s,\pi}[u^0] + \varepsilon_k \leq \frac{1}{2} \left( T \left\| Q \right\|_2 + \left\| Q_f \right\|_2 \right) |x_0|^2 e^{2C_A T} + \varepsilon_k .$$

On the other hand, since the matrix $R > 0$, one has

$$J_{s,\pi}[u^k] = \mathbb{E}_\pi \left[ \frac{1}{2} \int_s^T x_A^k(t)^T Q x_A^k(t) + u^k(t)^T R u^k(t) \, dt + \frac{1}{2} x_A^k(T)^T Q_f x_A^k(T) \right]$$
$$\geq \frac{1}{2} \int_s^T u^k(t)^T R u^k(t) \, dt \geq \frac{1}{2} \int_s^T r_1 |u^k(t)|^2 \, dt ,$$

where $r_1$ is the smallest eigenvalue of the matrix $R$.

Hence, one obtains the bound on the minimizing sequence

$$\int_s^T |u^k(t)|^2 \, dt \leq \frac{1}{r_1} \left( T \left\| Q \right\|_2 + \left\| Q_f \right\|_2 \right) |x_0|^2 e^{2C_A T} + \frac{2\varepsilon_k}{r_1} , \tag{24}$$

which results in a uniformly bounded norm:

$$\left\| u^k \right\|_2 \leq \sqrt{\frac{1}{r_1} \left( T \left\| Q \right\|_2 + \left\| Q_f \right\|_2 \right) |x_0|^2 e^{2C_A T} + \frac{2}{r_1}} =: C_u \quad \forall k \in \mathbb{N} . \tag{25}$$

In view of the previous relation, it follows from standard compactness arguments that $u^k \rightharpoonup \bar{u}$ weakly in $L^2([s, T]; \mathbb{R}^m)$. Since $u^k$ is uniformly bounded in $L^2$, then using in turn the relation (22), the Hölder inequality and the relation (25), one obtains that there exists a constant $C_x > 0$ such that

$$|x_A^k(t)| \leq C_x^{u^k} \leq \left( |x_0| + \sqrt{T} \left\| B \right\|_2 C_u \right) e^{C_A T} =: C_x \tag{26}$$

holds for every $k \in \mathbb{N}$, $A \in \mathcal{A}$ and $t \in [s, T]$. Furthermore, for each $k \in \mathbb{N}$, $A \in \mathcal{A}$ and $t \in [s, T]$, one has

$$\int_s^t \left| \dot{x}_A^k(\tau) \right| d\tau \leq \int_s^t \left\| A \right\|_2 |x_A^k(\tau)| d\tau + \int_s^t \left\| B \right\|_2 |u^k(\tau)| d\tau$$
$$\leq T C_A C_x + \sqrt{T} \left\| B \right\|_2 C_u ,$$

which implies that, for each $A \in \mathcal{A}$, $\dot{x}_A^k \rightharpoonup \dot{\bar{x}}_A$ weakly in $L^1([s, T]; \mathbb{R}^m)$, $x_A^k \to \bar{x}_A$ uniformly in $[s, T]$ and, in view of the linearity of the control system, the process $\{(\bar{x}_A, \bar{u})(\cdot) : A \in \mathcal{A}\}$ is the solution of the linear system

$$\begin{cases} \dot{\bar{x}}_A(t) = A\bar{x}_A(t) + B\bar{u}(t), & t \in [s, T], \\ \bar{x}_A(s) = x_0, \end{cases}$$

for each $A \in \mathcal{A}$. So the process $\{(\bar{x}_A, \bar{u})(\cdot) : A \in \mathcal{A}\}$ is a minimizer for Problem B, and in view of (24), $\bar{u}$ satisfies the bound (23) with the stricter constant

$$\overline{C}_u := \sqrt{\frac{1}{r_1} \left( T \, \|Q\|_2 + \|Q_f\|_2 \right) |x_0|^2 e^{2C_A T}} \,.$$

Since the functional $u \mapsto J_{s,\pi}[u]$ is strictly convex, the uniqueness of the minimizer follows from standard arguments. This completes the proof.                                    $\square$

**Remark 3.7** In view of the previous results, if $\{(\bar{x}_A, \bar{u})(\cdot) : A \in \mathcal{A}\}$ is an optimal process for problem B, then the constant

$$\overline{C}_x := \left( |x_0| + \sqrt{T} \, \|B\|_2 \, \overline{C}_u \right) e^{C_A T}$$

is such that

$$|\bar{x}_A(t)| \le \overline{C}_x, \quad \forall A \in \mathcal{A}, \quad \forall t \in [s, T]$$

and does not depend on $\pi$.

## 4 Optimality conditions

Let us consider the optimal control problem

$$\begin{cases} \text{minimize } J_{s,\pi}[u] \\ \text{over } u : [s, T] \to \mathbb{R}^m \text{ measurable such that} \\ u(t) \in U(t) \text{ a.e. } t \in [s, T], \\ \dot{x}_A(t) = A x_A(t) + B u(t), \quad A \in \mathcal{A}, \quad t \in [s, T] \\ x_A(s) = x_0, \end{cases} \tag{27}$$

where $U : [s, T] \rightsquigarrow \mathbb{R}^m$ is a $\mathcal{L} \times \mathcal{B}_{\mathbb{R}^m}$-measurable multifunction taking values compact sets and $J_{s,\pi}$ is the cost functional defined in (16) for a given, fixed $\pi \in \mathcal{M}(\mathcal{A})$. For a given reference process $\{(\bar{x}_A, \bar{u})(\cdot) : A \in \mathcal{A}\}$, we assume that the following condition holds true:

**(TH)**:  There exist $\delta > 0$ and a function $c \in L^2([s, T]; \mathbb{R})$ such that

$$|Ax + Bu| \le c(t),$$

for all $x \in \mathbb{B}_n(\bar{x}_A(t), \delta), u \in U(t), A \in \mathcal{A}$, a.e. $t \in [s, T]$.

It follows from standard ODE theory that when condition **(TH)** holds, for every admissible process $\{(\bar{x}_A, \bar{u})(\cdot) : A \in \mathcal{A}\}$ one has that $\bar{x}_A(\cdot)$ is in $W^{1,1}([s, T]; \mathbb{R}^n)$ for all $A \in \mathcal{A}$.

**Definition 4.1** For a given $\delta > 0$, a process $\{(\bar{x}_A, \bar{u})(\cdot) : A \in \mathcal{A}\}$ is said to be a $W^{1,1}$-local minimizer for problem (27) if

$$J_{s,\pi}[\bar{u}] \leq J_{s,\pi}[u]$$

for every process $\{(x_A, u)(\cdot) : A \in \mathcal{A}\}$ such that

$$\sup_{A \in \mathcal{A}} \|\bar{x}_A(\cdot) - x_A(\cdot)\|_{W^{1,1}} \leq \delta.$$

We recall the following result due to Bettiol and Khalil, which is a special case of Theorem 3.3 in [2]:

**Theorem 4.2** (Bettiol-Khalil, 2019) *Let $\{(\bar{x}_A, \bar{u})(\cdot) : A \in \mathcal{A}\}$ be a $W^{1,1}$-local minimizer for the optimal control problem* (27). *Let the assumption (TH) be satisfied. Then, there exists a $\mathcal{L} \times \mathcal{B}_{\mathcal{A}}$ measurable function $p : [s, T] \times \mathcal{A} \to \mathbb{R}^n$, $p(t, A) \equiv p_A(t)$, such that*

*(i)*

$$p_A(\cdot) \in W^{1,1}([s, T]; \mathbb{R}^n) \quad \forall A \in supp(\pi);$$

*(ii)*

$$\int_{\mathcal{A}} p_A(t) B \bar{u}(t)\, d\pi(A) - \frac{1}{2} \bar{u}(t)^T R \bar{u}(t)$$
$$= \max_{u \in U(t)} \left\{ \int_{\mathcal{A}} p_A(t) B u\, d\pi(A) - \frac{1}{2} u^T R u \right\} \quad \text{a.e. } t \in [s, T]$$

*(iii)*

$$-\dot{p}_A(t) = A^T p_A(t) - Q \bar{x}_A(t) \quad a.e.\, t \in [s, T],\, \forall A \in supp(\pi);$$

*(iv)*

$$-p_A(T) = Q_f \bar{x}_A(T) \quad \forall A \in supp(\pi).$$

**Remark 4.3** Let us notice that Theorem 3.3 in [2] is derived under the stronger assumption:

**(TH')**: there exist $\delta > 0$ and $c > 0$ such that

$$|Ax + Bu| \leq c,$$

for all $x \in \mathbb{B}(\bar{x}_A(t), \delta), u \in U(t), A \in \mathcal{A}$, a.e. $t \in [s, T]$.

However, scrutiny to the proof given there reveals that the result still holds true under the relaxed condition **(TH)**. Furthermore, Theorem 3.3 in [2] is derived for a Mayer optimal control problem, i.e., with only a final cost, but an analogous theorem for Bolza optimal control problems can be easily obtained by using a standard state augmentation argument.

We are now ready to prove the necessary optimality conditions for Problem B:

**Theorem 4.4** *Assume the hypothesis* **(SH)**. *Then, the following optimality condition is satisfied by the minimizer* $\{(\bar{x}_A, \bar{u})(\cdot) : A \in \mathcal{A}\}$ *for Problem B* (15). *There exists a* $\mathcal{L} \times \mathcal{B}_{\mathcal{A}}$ *measurable function* $p : [s, T] \times \mathcal{A} \to \mathbb{R}^n$, $p(t, A) \equiv p_A(t)$, *such that*

*(i)*

$$p_A(\cdot) \in W^{1,1}([s, T]; \mathbb{R}^n) \quad \forall A \in supp(\pi);$$

*(ii)*

$$\bar{u}(t) = +R^{-1}B^T \int_{\mathcal{A}} p_A(t) \, d\pi(A) \quad t \in [s, T];$$

*(iii)*

$$-\dot{p}_A(t) = A^T p_A(t) - Q\bar{x}_A(t) \quad a.e. \, t \in [s, T], \, \forall A \in supp(\pi);$$

*(iv)*

$$-p_A(T) = Q_f \bar{x}_A(T) \quad \forall A \in supp(\pi).$$

**Proof** Let us first observe that the optimal process $\{(\bar{x}_A, \bar{u})(\cdot) : A \in \mathcal{A}\}$ exists and is unique, due to Lemma 3.6. Consider now the optimal control problem

$$\begin{cases} \text{minimize } J_{s,\pi}[u] \\ \text{over } u : [s, T] \to \mathbb{R}^m \text{ measurable such that} \\ u(t) \in \mathbb{B}_m(\bar{u}(t), 1) \text{ a.e. } t \in [s, T], \\ \dot{x}_A(t) = Ax_A(t) + Bu(t), \quad A \in \mathcal{A}, \, t \in [s, T] \\ x_A(s) = x_0. \end{cases} \quad (28)$$

Such an optimal control problem is a special case of (27) with the choice of $U(t) = \mathbb{B}_m(\bar{u}(t), 1)$. Clearly, since $\{(\bar{x}_A, \bar{u})(\cdot) : A \in \mathcal{A}\}$ is a minimizer for Problem B, then it is also a minimizer for problem (28). Furthermore, since any element of $u \in U(t)$ can be written as $u = \bar{u}(t) + v$, for some $v \in \mathbb{B}_m(0, 1)$ and in view of Remark 3.7, then one can easily find $\delta > 0$ and a function $c \in L^2([s, T], \mathbb{R})$ such that the hypothesis **(TH)** is satisfied. Indeed,

$$|Ax + Bu| = |Ax + B(\bar{u}(t) + v)| \leq \|A\|_2 \, (\overline{C}_x + \delta) + \|B\|_2 \, (|\bar{u}(t)| + 1) =: c(t) \quad (29)$$

for all $x \in \mathbb{B}_n(\bar{x}_A(t), \delta)$, $\forall A \in \mathcal{A}$, for all $u \in \mathbb{B}_m(\bar{u}(t), 1)$, a.e. $t \in [s, T]$, where $\overline{C}_x$ is the constant appearing in Remark 3.7. So the process $\{(\bar{x}_A, \bar{u})(\cdot) : A \in \mathcal{A}\}$ is a $W^{1,1}$-minimizer (see Definition 4.1) for the optimal control problem (28), and the hypothesis **(SH)** is satisfied. Then, one can invoke Theorem 4.2, which provides conditions $(i)$–$(iii)$–$(iv)$ of the statement. In order to obtain condition $(ii)$, it is enough to observe that

$$\bar{u}(t) = \arg \max_{u \in U(t)} \left\{ \int_{\mathcal{A}} p_A(t) B u \, d\pi(A) - \frac{1}{2} u^T R u \right\}, \quad \text{a.e. } t \in [s, T]$$

and that $\bar{u}(t)$ is clearly an interior point of $U(t)$ for a.e. $t \in [s, T]$. Hence, $\bar{u}(t)$ has to satisfy also condition $(ii)$ of Theorem 4.4. This completes the proof.                    □

**Remark 4.5** Theorem 4.4 provides the existence of a multiplier $p_A(\cdot)$ for all $A \in supp(\pi)$. In general, we can extend its definition to all $A \in \mathcal{A}$, considering $p_A(\cdot)$ as the unique solution of

$$\begin{cases} -\dot{p}_A(t) = A^T p_A(t) - Q \bar{x}_A(t) & t \in [s, T] \\ -p_A(T) = Q_f \bar{x}_A(T), \end{cases} \quad (30)$$

where $\{(\bar{x}_A, \bar{u})(\cdot) : A \in \mathcal{A}\}$ is the unique minimizer of Problem B.

The following result guarantees that this extension defined in Remark 4.5 is continuous with respect to $A$:

**Lemma 4.6** (Boundness and continuity of multipliers) *Let us consider the multipliers defined in* (30) *for each* $A \in \mathcal{A}$. *They have the following properties:*

*(i) There exists a positive constant* $\overline{C}_p$, *independent from* $\pi$, *which bounds uniformly all multipliers, i.e.,*

$$|p_A(t)| \le \overline{C}_p \quad \forall t \in [s, T], \forall A \in \mathcal{A},$$

*(ii) The map* $A \mapsto p_A(\cdot)$ *is continuous.*

**Proof** The proof is similar to that of Lemma 3.5, with the only difference that here we need to apply the Grönwall lemma backward instead of forward. Notice that the final condition

$$-p_A(T) = Q_f \bar{x}_A(T)$$

will not be the same for all $A \in \mathcal{A}$, whereas the initial condition was the same in Lemma 3.5. However, it is still continuous with respect to $A$ by Lemma 3.5, so the same arguments can be applied to prove $(ii)$.                    □

## 5 Main convergence results

In this section, we will present the main results, which are, respectively, the convergence of the value function (Corollary 5.2) and the convergence of the optimal control (Theorem 5.3). Given a sequence of probability distributions $\{\pi^N\} \subset \mathcal{M}(\mathcal{A})$, for each $N \in \mathbb{N}$ we consider *Problem $B_{\pi^N}$*, namely problem (15) relative to the distribution $\pi^N$. Recalling the definition of the value function in (17), we define

$$V^N(s, x_0) := V_{\pi^N}(s, x_0) = \inf_{u \in \mathcal{U}} J_{s, \pi^N}[u] \,. \tag{31}$$

If $\{\pi^N\} \subset \mathcal{M}(\mathcal{A})$ is such that $W_1(\pi^N, \pi^\infty) \to 0$ for $N \to \infty$, what can be said about the convergence of the value functions $V^N$ to $V^\infty$ and of the optimal controls $u^N$ to $u^\infty$ for $N \to \infty$? In this section, we give an answer to those questions.

**Theorem 5.1** (Lipschitz estimate for the value function w.r.t. $\pi$) *Let the assumption (SH) be satisfied. Given $\pi, \pi' \in \mathcal{M}(\mathcal{A})$ and $(s, x_0) \in [0, T] \times K$ with $K \subset \mathbb{R}^n$ compact set, let us consider the two value functions $V_\pi$ and $V_{\pi'}$ as defined in (17). Then, the distance between $V$ and $V'$ can be bounded uniformly for $(s, x_0) \in [0, T] \times K$, that is:*

$$\|V_\pi - V_{\pi'}\|_{\infty, [0,T] \times K} \le C_K \, W_1(\pi, \pi') \,, \tag{32}$$

*where $C_K = C_K(T, C_A, \|Q\|_2, \|Q_f\|_2, r_1, K)$ is a constant which does not depend on the distributions $\pi$ and $\pi'$, but merely on the compact set $\mathcal{A}$.*

**Proof** We divide the proof into three steps.

*STEP 1*: Fix two matrices $A, A' \in \mathcal{A}$, any point $(s, x_0) \in [0, T] \times \mathbb{R}^n$ and a control $u \in \mathcal{U}$. Using Grönwall's lemma as we did for point (ii) of Lemma 3.5, we get the following estimate:

$$\begin{aligned}
|x_A(t) - x_{A'}(t)| &\le C_x^u \, (t - s) \, e^{\|A\|_2 (t-s)} \, \|A - A'\|_2 \\
&\le C_x^u \, t \, e^{\|A\|_2 t} \, \|A - A'\|_2 \qquad \forall \, t \in [s, T] \,,
\end{aligned}$$

with $C_x^u$ given by Lemma 3.5.

Let us denote by $\ell$ the running cost and by $h$ the final cost:

$$\ell(x, u) := \frac{1}{2} \left( x^T Q x + u^T R u \right) \quad \text{and} \quad h(x) := \frac{1}{2} x^T Q_f x \,,$$

so we can write

$$J_{s, \pi}[u] = \int_{\mathcal{A}} \left[ \int_s^T \ell(x_A(t), u(t)) \, dt + h(x_A(T)) \right] d\pi(A) \,.$$

Since both $\ell$ and $h$ are locally Lipschitz continuous, one has

$$
\begin{aligned}
|\ell(x_{A'}(t), u(t)) - \ell(x_A(t), u(t))| &= \frac{1}{2}|x_{A'}{}^T(t) Q x_{A'}(t) - x_A^T(t) Q x_A(t)| \\
&\leq \frac{1}{2}|x_{A'}{}^T Q(x_{A'}(t) - x_A(t))| \\
&\quad + \frac{1}{2}|(x_{A'}(t) - x_A(t))^T Q x_A(t)| \\
&\leq \|Q\|_2 \, C_x |(x_{A'}(t) - x(t))| \\
&\leq L_\ell^u \, C_x^u \, t \, e^{\|A\|_2 t} \, \left\|A - A'\right\|_2 ,
\end{aligned}
$$

and, similarly,

$$
|h(x_{A'}(T)) - h(x_A(T))| \leq L_h^u \, C_x^u \, T \, e^{\|A\|_2 T} \, \left\|A - A'\right\|_2 ,
$$

where we defined the Lipschitz constants

$$
L_\ell^u := \|Q\|_2 \, C_x^u \quad \text{and} \quad L_h^u := \left\|Q_f\right\|_2 \, C_x^u ;
$$

these two constants inherit from $C_x^u$ the dependency on $x_0$ and $u$.

Finally, the cost difference between two single trajectories can be easily bounded by

$$
\begin{aligned}
\int_s^T &\left|\ell(x_{A'}(t), u(t)) - \ell(x_A(t), u(t))\right| dt + |h(x_{A'}(T)) - h(x_A(T))| \\
&\leq \left(T L_\ell^u + L_h^u\right) C_x^u T \, e^{C_A T} \, \left\|A - A'\right\|_2 .
\end{aligned}
\tag{33}
$$

*STEP 2*: Fix an initial condition $x(s) = x_0 \in \mathbb{R}^n$ with $s \in [0, T]$ and a control $u \in \mathcal{U}$. We want to prove a bound for the distance between $J_{s,\pi}[u]$ and $J_{s,\pi'}[u]$.

As a property of $W_1$, there exists (see Theorem 4.1 on [27]) a probability distribution $\gamma^* \in \Gamma(\pi, \pi')$ on $\mathcal{A} \times \mathcal{A}$ with marginal distributions $\pi$ and $\pi'$ such that

$$
W_1(\pi, \pi') = \int_{\mathcal{A} \times \mathcal{A}} d_2(A, A') \, d\gamma^*(A, A') ,
\tag{34}
$$

where $d_2$ is the distance introduced in (7). We can thus write

$$
\begin{aligned}
\left| J_{s,\pi}[u] - J_{s,\pi'}[u] \right| &= \left| \int_{\mathcal{A}} \int_s^T \ell(x_A(t), u(t)) \, dt \, d\pi(A) \right. \\
&\quad - \int_{\mathcal{A}} \int_s^T \ell(x_{A'}(t), u(t)) dt d\pi'(A') \\
&\quad \left. + \int_{\mathcal{A}} h(x_A(T)) d\pi(A) - \int_{\mathcal{A}} h(x_{A'}(T)) d\pi'(A') \right| \\
&= \left| \int_{\mathcal{A} \times \mathcal{A}} \left[ \int_s^T \left( \ell(x_A(t), u(t)) - \ell(x_{A'}(t), u(t)) \right) \, dt \right. \right. \\
&\quad \left. \left. + h(x_A(T)) - h(x_{A'}(T)) \right] d\gamma^*(A, A') \right|,
\end{aligned}
$$

where we have used that

$$
\int_{\mathcal{A} \times \mathcal{A}} \left( \phi(A) + \psi(A') \right) d\gamma^*(A, A') = \int_{\mathcal{A}} \phi(A) d\pi(A) + \int_{\mathcal{A}} \psi(A') d\pi'(A')
$$

for all measurable functions $\phi, \psi$ on $\mathcal{A}$, since $\gamma^*$ admits $\pi$ and $\pi'$ as marginals.

Using the bound (33) from STEP 1 and formula (34), we get

$$
\begin{aligned}
\left| J_{s,\pi}[u] - J_{s,\pi'}[u] \right| &\le \int_{\mathcal{A} \times \mathcal{A}} \left( T L_\ell^u + L_h^u \right) C_x^u \, T \, e^{C_A T} \left\| A - A' \right\|_2 \, d\gamma^*(A, A') \\
&= \left( T L_\ell^u + L_h^u \right) C_x^u \, T \, e^{C_A T} \, W_1(\pi, \pi') .
\end{aligned} \tag{35}
$$

Note that the constant $C_x^u$ which appears here depends merely on $x_0$ and $u$.

*STEP 3*: We will now show that an estimate similar to (35) holds true even for the value functions $V_\pi$ and $V_{\pi'}$.

Fix any point $(s, x_0) \in [0, T] \times K$. In view of Lemma 3.6, there exist controls $\bar{u}, \bar{u}' \in \mathcal{U}$ such that

$$
J_{s,\pi}[\bar{u}] = V_\pi(s, x_0), \qquad J_{s,\pi'}[\bar{u}'] = V_{\pi'}(s, x_0) .
$$

Then, one has

$$
V_{\pi'}(s, x_0) - V_\pi(s, x_0) = \inf_{u \in \mathcal{U}} J_{s,\pi'}[u] - J_{s,\pi}[\bar{u}] \le J_{s,\pi'}[\bar{u}] - J_{s,\pi}[\bar{u}]
$$

and, in the same way,

$$
V_\pi(s, x_0) - V_{\pi'}(s, x_0) = \inf_{u \in \mathcal{U}} J_{s,\pi}[u] - J_{s,\pi'}[\bar{u}] \le J_{s,\pi}[\bar{u}'] - J_{s,\pi'}[\bar{u}'] .
$$

Hence,

$$|V_\pi(s, x_0) - V_{\pi'}(s, x_0)| \le \max \left\{ \left| J_{s,\pi'}[\bar{u}] - J_{s,\pi}[\bar{u}] \right|, \left| J_{s,\pi'}[\bar{u}'] - J_{s,\pi}[\bar{u}'] \right| \right\}.$$

Moreover, being both $\bar{u}$ and $\bar{u}'$ optimal control for some distribution on $\mathcal{A}$, we can use the uniform constant $\overline{C}_x$ given by Remark 3.7 and define an analogous uniform constants for the Lipschitz continuity of $\ell$ and $h$:

$$\overline{L}_\ell := \|Q\|_2 \, \overline{C}_x, \quad \overline{L}_h := \|Q_f\|_2 \, \overline{C}_x.$$

In this way, the estimate becomes independent from $\bar{u}$ and $\bar{u}'$ and thus from $\pi$ and $\pi'$:

$$|V_{\pi'}(s, x_0) - V_\pi(s, x_0)| \le C_{x_0} \, W_1(\pi, \pi'),$$

where

$$C_{x_0} := \left( T \, \overline{L}_\ell + \overline{L}_h \right) \overline{C}_x T \, e^{C_A T}.$$

Finally, noting that the estimate depends on $x_0$ only through its norm, we get that the bound is uniform in compact sets $K \subset \mathbb{R}^n$, letting

$$C_K := \sup_{x_0 \in K} C_{x_0}.$$

$\square$

A straightforward consequence of the previous theorem is the following:

**Corollary 5.2** (Convergence of the value function) *If* $\pi^N \overset{*}{\rightharpoonup} \pi^\infty$, *then* $V^N(s, x_0) \to V^\infty(s, x_0)$ *for each* $(s, x_0) \in [0, T] \times \mathbb{R}^n$. *The convergence is uniform in compact sets* $[0, T] \times K$, *where* $K \subset \mathbb{R}^n$ *is compact.*

In what follows, we will use $\bar{u}^N(\cdot)$ to denote the optimal control of *Problem* $B_{\pi^N}$. Furthermore, $\bar{x}_A^N(\cdot)$ and $p_A^N(\cdot)$ denote, respectively, the optimal trajectories and the multipliers relative to *Problem* $B_{\pi^N}$, that is the solutions of the differential systems (18) and (30).

The following theorem provides a strong convergence of $\bar{u}^N(\cdot)$ to the optimal control $\bar{u}^\infty(\cdot)$ of the limit problem $B_{\pi^\infty}$, assuming that $\pi^N \overset{*}{\rightharpoonup} \pi^\infty$.

**Theorem 5.3** (Convergence of the optimal control) *Let the assumption* **(SH)** *be satisfied. Consider a sequence of probability distributions* $\{\pi^N\} \subset \mathcal{M}(\mathcal{A})$, *such that* $\pi^N \overset{*}{\rightharpoonup} \pi^\infty$ *and fix* $s \in [0, T]$ *and* $x_0 \in \mathbb{R}^n$. *If* $\bar{u}^N(\cdot)$ *and* $\bar{u}^\infty(\cdot)$ *are, respectively, the optimal controls of Problem* $B_{\pi^N}$ *and* $B_{\pi^\infty}$, *i.e., they satisfy* (19), *respectively, for* $\pi^N$ *and* $\pi^\infty$, *then* $\bar{u}^N(\cdot)$ *converges uniformly* $\bar{u}^\infty(\cdot)$ *for* $N \to \infty$.

*Proof* Without loss of generality, one can take $s = 0$, being all the other cases similar.

Lemma 3.6 assures that, for each $\pi^N \in \mathcal{M}(\mathcal{A})$, there exists a unique optimal process $\left\{ (\bar{x}_A^N, \bar{u}^N)(\cdot) : A \in \mathcal{A} \right\}$. Taking into account Theorem 4.4 and Remark 4.5, that process satisfies the following necessary conditions: For each $N \in \mathbb{N}$, there exists a continuous function $p^N : [0, T] \times \mathcal{A} \to \mathbb{R}^n$, $p^N(t, A) \equiv p_A^N(t)$, such that

(i)

$$p_A^N(\cdot) \in W^{1,1}([0,T]; \mathbb{R}^n) \quad \forall A \in \mathcal{A};$$

(ii)

$$\bar{u}^N(t) = +R^{-1} B^T \int_{\mathcal{A}} p_A^N(t) \, d\pi^N(A), \quad \forall \, t \in [0,T];$$

(iii)

$$-\dot{p}_A^N(t) = A^T p_A^N(t) - Q \bar{x}_A^N(t) \quad a.e. \ t \in [0,T], \ \forall A \in \mathcal{A};$$

(iv)

$$-p_A^N(T) = Q_f \bar{x}_A^N(T) \quad \forall A \in \mathcal{A}.$$

For each $A \in \mathcal{A}$, the family of functions

$$\mathcal{F}_A := \left\{ \left( \bar{x}_A^N, p_A^N \right)(\cdot) \right\}_{N \in \mathbb{N}}$$

is uniformly bounded due to Lemmas 3.5 and 4.6. Moreover, one can find the bounds on the derivatives:

$$\int_0^T |\dot{\bar{x}}_A^N(t)| dt = \int_0^T \left| A \bar{x}_A^N(t) + B \bar{u}^N(t) \right| dt \leq T C_A \overline{C}_x + \sqrt{T} \, \|B\|_2 \, \overline{C}_u, \tag{36}$$
$$|\dot{p}_A^N(t)| = |A^T p_A^N(t) - Q \bar{x}_A^N(t)| \leq C_A \overline{C}_p + \|Q\|_2 \, \overline{C}_x,$$

for all $A \in \mathcal{A}$, a.e. $t \in [0,T]$. The second bound in (36) and the relation (ii) imply that also the map

$$[0,T] \ni t \mapsto \bar{u}^N(t) = +R^{-1} B^T \int_{\mathcal{A}} p_A^N(t) \, d\pi^N(A)$$

is equibounded and equicontinuous in $N$. For each $A \in \mathcal{A}$ fixed, one can then apply Theorem 2.5.3 on [28], implying the existence of some limit functions $x_A^\infty, p_A^\infty \in W^{1,1}([0,T]; \mathbb{R}^n)$ and $u^\infty \in C^0([0,T]; \mathbb{R}^m)$ such that

$$\dot{\bar{x}}_A^N(t) \rightharpoonup \dot{x}_A^\infty(t) \text{ and } \dot{p}_A^N(t) \rightharpoonup \dot{p}_A^\infty(t) \quad \text{weakly in } L^1([0,T]; \mathbb{R}^n) \text{ as } N \to \infty,$$
$$\bar{x}_A^N(t) \to x_A^\infty(t) \quad \text{and} \quad p_A^N(t) \to p_A^\infty(t) \quad \text{uniformly on } [0,T] \text{ as } N \to \infty,$$
$$\bar{u}^N(t) \to u^\infty(t) \quad \text{uniformly on } [0,T] \text{ as } N \to \infty,$$

and such that, for each $A \in \mathcal{A}$, $(x_A^\infty, p_A^\infty)(\cdot)$ is a solution of the boundary value problem

$$
\begin{cases}
\dot{x}_A^\infty(t) = A x_A^\infty(t) + B u^\infty(t), & t \in [0, T] \\
-\dot{p}_A^\infty(t) = A^T p_A^\infty(t) - Q x_A^\infty(t), & t \in [0, T] \\
x_A^\infty(0) = x_0 \\
-p_A^\infty(T) = Q_f x_A^\infty(T).
\end{cases}
\tag{37}
$$

Furthermore, since $\pi^N \overset{*}{\rightharpoonup} \pi^\infty$, $u^\infty$ satisfies the relation

$$
[0, T] \ni t \mapsto u^\infty(t) = +R^{-1}B^T \int_{\mathcal{A}} p_A^\infty(t)\, d\pi^\infty(A).
\tag{38}
$$

In fact, the result follows from the estimate

$$
\left| \int_{\mathcal{A}} p_A^N(t)\, d\pi^N(A) - \int_{\mathcal{A}} p_A^\infty(t)\, d\pi^\infty(A) \right| \le \left| \int_{\mathcal{A}} p_A^N(t)\, d\pi^N(A) - \int_{\mathcal{A}} p_A^N(t)\, d\pi^\infty(A) \right|
$$
$$
+ \left| \int_{\mathcal{A}} p_A^N(t)\, d\pi^\infty(A) - \int_{\mathcal{A}} p_A^\infty(t)\, d\pi^\infty(A) \right| \quad (39)
$$

for all $t \in [0, T]$, which implies that

$$
\bar{u}^N(t) = R^{-1}B^T \int_{\mathcal{A}} p_A^N(t)\, d\pi^N(A) \to u^\infty(t) = R^{-1}B^T \int_{\mathcal{A}} p_A^\infty(t)\, d\pi^\infty(A)
$$

uniformly on $[0, T]$ as $N \to \infty$.

Notice that the convergence is guaranteed along a subsequence, but one can say that the whole sequence converges since the limit does not depend on the subsequence. (It solves (37).)

It remains to show that the limiting process $\{(x_A^\infty, u^\infty)(\cdot) : A \in \mathcal{A}\}$ is actually optimal for the Problem B (15) relative to $\pi^\infty$. To this aim, let us stress the following properties of the cost functional $J_{0,\pi}$ in (16), using the lighter notation

$$
J^N := J_{0,\pi^N} \quad \text{and} \quad J^\infty := J_{0,\pi^\infty} :
$$

1) if $\pi^N \overset{*}{\rightharpoonup} \pi^\infty$, $J^N[u] \to J^\infty[u]$ for each $u \in \mathcal{U}$, since the map $A \mapsto x_A(\cdot)$ is continuous by Lemma 3.5;
2) each $J^N$ is continuous with respect to $u$, since for each $A \in \mathcal{A}$, the map $u \mapsto x_A(\cdot; u)$ is continuous, again from Lemma 3.5.

Since $\bar{u}^N$ is the optimal control of Problem $B_{\pi^N}$ and $u$ is an admissible control for the same problem, then we get

$$
J^N[\bar{u}^N] \le J^N[u], \quad \forall N \in \mathbb{N}
$$

so, letting $N \to \infty$, it easily follows from the previous relation and properties 1) and 2) that

$$J^\infty[u^\infty] \le J^\infty[u] .$$

In view of the uniqueness of the optimal control $\bar{u}^\infty$ for Problem $B_{\pi^\infty}$, one can conclude that $u^\infty \equiv \bar{u}^\infty$. Hence, also the process $\left\{ (x_A^\infty, u^\infty)(\cdot) : A \in \mathcal{A} \right\}$ is optimal for the given problem. This concludes the proof. $\quad\square$

**Remark 5.4** *(Special Case: $\pi = \delta_{\hat{A}}$)* The particular case in which $\pi$ is a Dirac delta $\delta_{\hat{A}}$ for a given matrix $\hat{A} \in \mathbb{M}_{n \times n}$ deserves special attention. Indeed, when $\pi = \delta_{\hat{A}}$, the cost functional $J_{s,\pi}$ (16) becomes the cost functional

$$J_s[u] := \frac{1}{2} \int_s^T \left( x(t)^T Q x(t) + u(t)^T R u(t) \right) dt + \frac{1}{2} x(T)^T Q_f x(T) ,$$

and Problem B in (15) coincides with a standard Problem A (see (10))

$$\begin{cases} \text{minimize } \left\{ \frac{1}{2} \int_s^T \left( x(t)^T Q x(t) + u(t)^T R u(t) \right) dt + \frac{1}{2} x(T)^T Q_f x(T) \right\} \\ \text{over } (x, u)(\cdot) \text{ such that } u : [s, T] \to \mathbb{R}^m \text{ is measurable and} \\ \dot{x}(t) = \hat{A} x(t) + B u(t), \quad t \in [s, T] \\ x(s) = x_0 . \end{cases} \quad (40)$$

Furthermore, the definition of the value function and that of the optimal control we gave in (17) and (19) agree, in this special case, with the classic definitions in control theory (see (13) and (14)):

$$V(s, x_0) := \inf_{u \in \mathcal{U}} J_s[u] \text{ and } \bar{u} := \arg\min_{u \in \mathcal{U}} J_s[u] .$$

If we apply Corollary 5.2 and Theorem 5.3 to a sequence $\pi^N$ converging to $\delta_{\hat{A}}$, then we obtain

$$V^N(s, x_0) \to V(s, x_0) \quad \forall s \in [0, T], \ x_0 \in \mathbb{R}^n$$

and

$$\bar{u}^N \to \bar{u} \quad \text{uniformly in } [0, T] ,$$

where $V^N$ and $\bar{u}^N$ are, respectively, the value function in (17) and the optimal control (19) relative to $\pi^N$.

**Remark 5.5** *(Additional remarks on the connection with RL)* Let us comment the results of this section bearing in mind the PILCO algorithm, presented in Remark 3.4. Recall that PILCO uses the agent experience on the environment to tune the parameters of a Gaussian process (GP), building up a stochastic model for the dynamics. The GP is then updated as the agent collects more data on the environment. This procedure generates a sequence of probability distributions $\pi^N$ on the space of continuous functions, which should get sufficiently close to the Dirac delta representing the actual dynamics.

Now, assume that we are applying PILCO to the linear system (12) in which the agent does not know the matrix $\hat{A}$. Then, one can think to $\left\{\pi^N\right\}_{N\in\mathbb{N}}$ as a sequence of probability distributions on the space of matrices $\mathbb{M}_{n\times n}$, which is converging to a Dirac delta concentrated at the true matrix $\hat{A}$. If, at each step, the agent picks as control $u^N$, namely the one minimizing the expected cost with respect to $\pi^N$, and $\pi^N \overset{*}{\rightharpoonup} \pi^\infty$, then one can apply Theorem 5.3 to the special case presented in Remark 5.4 and obtains the uniform convergence $\bar{u}^N \to \bar{u}^\infty$, where $\bar{u}^\infty$ is the optimal control for the limit problem.

This suggests that even if the distribution $\pi^N$ does not exactly reach $\delta_{\hat{A}}$ but gets sufficiently close to it, then the control $\bar{u}^N$ is suboptimal for the actual LQR problem.

## 6 A case of study: finite support measures converging to $\delta_{\hat{A}}$

In this section, we will assume that $\mathcal{A}$ is a finite set, namely $\mathcal{A} := \{A_1, \ldots, A_M\}$ for some integer $M \in \mathbb{N}$. Let us consider a sequence of probability distributions $\{\pi^N\} \subset \mathcal{M}(\mathcal{A})$, which can be written as

$$\pi^N := \sum_{i=1}^M \alpha_i^N \delta_{A_i}, \quad \text{for some } \alpha_i^N \geq 0 \text{ such that } \sum_{i=1}^M \alpha_i^N = 1, \quad \forall N \in \mathbb{N}. \quad (41)$$

For a given $s \in [0, T]$ and $x_0 \in \mathbb{R}^n$, suppose that the underlying dynamics governing the optimal control problem we are interested in is a standard Problem A (see (10)):

$$\begin{cases} \text{minimize } J_s[u] \\ \text{over } (x, u)(\cdot) \text{ such that } u : [s, T] \to \mathbb{R}^m \text{ is measurable and} \\ \dot{x}(t) = \hat{A}x(t) + Bu(t), \qquad t \in [s, T] \\ x(s) = x_0, \end{cases} \quad (42)$$

where $\hat{A} \in \mathcal{A}$ and the cost functional $J_s[u]$ is defined as in (11):

$$J_s[u] := \frac{1}{2} \int_s^T \left( x(t)^T Q x(t) + u(t)^T R u(t) \right) dt + \frac{1}{2} x(T)^T Q_f x(T). \quad (43)$$

Without loss of generality, one can set $\hat{A} \equiv A_1$. For each $s \in [0, T]$ and $x_0 \in \mathbb{R}^n$, the value function $V(s, x_0)$ for this problem has been defined in (13).

Then, one can expect that, after some interactions with the system, it is possible to construct a sequence of probability distributions $\{\pi^N\} \subset \mathcal{M}(\mathcal{A})$ capturing the current belief that one has about the real system (42) and such that, when $N$ is sufficiently large, $\pi^N$ gets closer and closer to $\delta_{A_1}$ and eventually $\pi^N \overset{*}{\rightharpoonup} \delta_{A_1}$.

For each fixed $N \in \mathbb{N}$, one can reformulate Problem B associated with $\pi^N$ as a classical LQR problem on an augmented system of dimension $nM$:

$$\begin{cases} \text{minimize } J_s^N[u] \\ \text{over } (X, u)(\cdot) \text{ such that } u : [s, T] \to \mathbb{R}^m \text{ is measurable and} \\ \dot{X}(t) = \widetilde{A} X(t) + \widetilde{B} u(t), \quad t \in [s, T] \\ X(s) = X_0 , \end{cases} \tag{44}$$

with cost functional

$$J_s^N[u] := \frac{1}{2} \int_s^T \left( X(t)^T \widetilde{Q}^N X(t) + u(t)^T R u(t) \right) dt + \frac{1}{2} X(T)^T \widetilde{Q}_f^N X(T) . \tag{45}$$

where we have used the compact notation $X(t) := \left( x_{A_1}(t), \ldots, x_{A_M}(t) \right)$, $X_0 \in \mathbb{R}^{nM}$ and

$$\widetilde{A} = \begin{pmatrix} A_1 & 0 & \cdots & 0 \\ 0 & A_2 & \cdots & 0 \\ \vdots & \vdots & \ddots & \vdots \\ 0 & 0 & \cdots & A_M \end{pmatrix}, \quad \widetilde{B} = \begin{pmatrix} B \\ B \\ \vdots \\ B \end{pmatrix},$$

$$\widetilde{Q}^N = \begin{pmatrix} \alpha_1^N Q & 0 & \cdots & 0 \\ 0 & \alpha_2^N Q & \cdots & 0 \\ \vdots & \vdots & \ddots & \vdots \\ 0 & 0 & \cdots & \alpha_M^N Q \end{pmatrix} \text{ and } \widetilde{Q}_f^N = \begin{pmatrix} \alpha_1^N Q_f & 0 & \cdots & 0 \\ 0 & \alpha_2^N Q_f & \cdots & 0 \\ \vdots & \vdots & \ddots & \vdots \\ 0 & 0 & \cdots & \alpha_M^N Q_f \end{pmatrix}.$$

In this section, we will use $V^N(s, X_0)$ to denote the value function related to the LQR problem (44), namely

$$V^N(s, X_0) = \inf_{u \in \mathcal{U}} J_s^N[u]. \tag{46}$$

Since the optimal control problem (44) can be regarded as a classic LQR problem, then one has the following relation between the value function and the optimal control in feedback form (see, e.g., [13], Theorem 3.4):
there exists $P^N$ such that

$$V^N(s, X_0) = X_0^T P^N(s) X_0, \quad \nabla_X V^N(s, X_0) = 2 P^N(s) X_0 \tag{47}$$

$$\bar{u}^N(s, X_0) = -R^{-1} \widetilde{B}^T P^N(s) X_0, \tag{48}$$

where $[s, T] \ni t \mapsto P^N(t) \in \mathbb{M}_{nM \times nM}$ solves the Riccati equation

$$\begin{cases} \widetilde{A}^T P(t) + P(t)\widetilde{A} - P(t)\widetilde{B}R^{-1}\widetilde{B}^T P(t) + \widetilde{Q}^N = -\dot{P}(t), & t \in [s, T] \\ P(T) = \widetilde{Q}_f^N. \end{cases}$$

$$(49)$$

For $x_0 \in \mathbb{R}^n$, we will use the notation $V^N(s, x_0)$ to denote the value function $V^N(s, X_0)$ evaluated at $X_0 = (x_0, \ldots, x_0) \in \mathbb{R}^{nM}$.

We can summarize the results of the previous section applied to problem (44) as follows:

**Corollary 6.1** *Let the assumption **(SH)** be satisfied. For each $s \in [0, T]$, $x_0 \in \mathbb{R}^n$ and $\{\pi^N\} \subset \mathcal{M}(\mathcal{A})$, the optimal control problems (44) and (42) satisfy the following conditions:*

*(i) problems (44) and (42) admit a unique optimal process $\left\{ (\bar{x}_A^N, \bar{u}^N)(\cdot) : A \in \mathcal{A} \right\}$ and $(\bar{x}, \bar{u})(\cdot)$, respectively;*

*(ii) for each $K \subset \mathbb{R}^n$, $\exists C_K$ such that*

$$||V^N - V||_{\infty, K} \le C_K W_1(\pi^N, \delta_{\hat{A}}) ;$$

$$(50)$$

*(iii) if, moreover, $\{\pi^N\}$ is such that $\pi^N \overset{*}{\rightharpoonup} \delta_{\hat{A}}$, then the optimal control $\bar{u}^N \to \bar{u}$ uniformly for $t \in [s, T]$.*

For each $N \in \mathbb{N}$, the solution of (49) is related to the matrices $X^N, Y^N : [s, T] \to \mathbb{M}_{nM \times nM}$, such that the pair $(X^N, Y^N)(\cdot)$ is the solution of the backward Hamiltonian differential equation

$$\begin{cases} \begin{bmatrix} \dot{X}(t) \\ \dot{Y}(t) \end{bmatrix} = H^N \begin{bmatrix} X(t) \\ Y(t) \end{bmatrix} & \text{for } t \in [s, T] \\ \begin{bmatrix} X(T) \\ Y(T) \end{bmatrix} = \begin{bmatrix} I \\ \widetilde{Q}_f^N \end{bmatrix}, \end{cases}$$

$$(51)$$

where

$$H^N = \begin{bmatrix} \widetilde{A} & -\widetilde{B}R^{-1}\widetilde{B}^T \\ \widetilde{Q}^N & -\widetilde{A}^T \end{bmatrix}.$$

$$(52)$$

The relation between the solutions of (49) and (51) was stated in precise terms by Coppel in [5, pp. 274–275]:

**Theorem 6.2** *Suppose that **(SH)** holds true. Let $X, Y : [s, T] \to \mathbb{M}_{nM \times nM}$ be the solutions of the Hamiltonian differential problem (51). Then,*

*1. $X(t)$ is non-singular for all $t \in [s, T]$;*
*2. the solution of (49) is*

$$\widetilde{P}(t) = Y(t)X^{-1}(t), \quad t \in [s, T].$$

$$(53)$$

We are now ready to strengthen the results of the previous by showing that, in the case in which the sequence of measures and its limit have finite support, then also the solution of the Riccati equation (49) has good convergence properties:

**Theorem 6.3** *(Convergence of the Riccati equation) Let the assumption **(SH)** be satisfied. Suppose that the sequence $\{\pi^N\} \subset \mathcal{M}(\mathcal{A})$ is such that $\pi^N \overset{*}{\rightharpoonup} \delta_{A_1}$, that is, for each $i = 1, \ldots, M$ the weights converge:*

$$\begin{aligned}
\alpha_1^N &\to 1 \\
\alpha_i^N &\to 0 \qquad for \; i = 2, \ldots, M.
\end{aligned}$$

*when $N \to \infty$. Then, the sequence of matrices $\{\widetilde{P}^N(t)\} \subset \mathbb{M}_{nM \times nM}$ which solve the Riccati equation (49) for each $N \in \mathbb{N}$ converges to the matrix*

$$\bar{P}(t) = \begin{pmatrix} P(t) & 0 & \cdots & 0 \\ 0 & 0 & \cdots & 0 \\ \vdots & \vdots & \ddots & \vdots \\ 0 & 0 & \cdots & 0 \end{pmatrix},$$

*uniformly on $t \in [s, T]$, where $P(t) \in \mathbb{M}_{n \times n}$ is the solution of the Riccati equation related to the optimal control problem (42) with state matrix $A_1$, namely:*

$$\begin{cases} A_1^T P(t) + P(t)A_1 - P(t)BR^{-1}B^T P(t) + Q = -\dot{P}(t), & t \in [s, T], \\ \qquad\qquad\qquad\qquad\qquad\qquad\qquad P(T) = Q_f\,, \end{cases} \tag{54}$$

**Proof** Consider the Hamiltonian systems in (51) related to all different distributions $\pi^N$. Notice that the norm of $H^N$ in (52) can be bounded:

$$\left\| H^N \right\|_2 \leq 2C_A + \|Q\|_2 + M \ \|B\|_2^2 \left\| R^{-1} \right\|_2,$$

for all $N \in \mathbb{N}$.

Using Grönwall's lemma, one can easily show that the pair of matrices $(X^N, Y^N)$ solution to (51) is uniformly bounded and that, using again (51), $(\dot{X}^N, \dot{Y}^N)$ is uniformly integrally bounded, for all $N \in \mathbb{N}$. So it is possible to apply Theorem 2.5.3 of [28] to show that the pair $(X^N, Y^N)$ converges to some matrices $(X^\infty, Y^\infty)$ solution of the system

$$\begin{cases} \begin{bmatrix} \dot{X}(t) \\ \dot{Y}(t) \end{bmatrix} = H^\infty \begin{bmatrix} X(t) \\ Y(t) \end{bmatrix} & \text{for } t \in [s, T] \\ \begin{bmatrix} X(T) \\ Y(T) \end{bmatrix} = \begin{bmatrix} I \\ \widetilde{Q}_f^\infty \end{bmatrix}, \end{cases} \tag{55}$$

where

$$H^\infty = \begin{bmatrix} \widetilde{A} & -\widetilde{B}R^{-1}\widetilde{B}^T \\ \widetilde{Q}^\infty & -\widetilde{A}^T \end{bmatrix} \tag{56}$$

and

$$\widetilde{Q}^\infty = \begin{pmatrix} Q & 0 & \cdots & 0 \\ 0 & 0 & \cdots & 0 \\ \vdots & \vdots & \ddots & \vdots \\ 0 & 0 & \cdots & 0 \end{pmatrix}, \quad \widetilde{Q}_f^\infty = \begin{pmatrix} Q_f & 0 & \cdots & 0 \\ 0 & 0 & \cdots & 0 \\ \vdots & \vdots & \ddots & \vdots \\ 0 & 0 & \cdots & 0 \end{pmatrix}.$$

Then, in view of Theorem 6.2, the matrix $X^\infty(t)$ is also nonsingular for each $t \in [s, T]$, and the matrix $\widetilde{P}^\infty(t) := Y^\infty(t)X^\infty(t)^{-1}$ is the solution of the Riccati equation

$$\begin{cases} \widetilde{A}^T\widetilde{P}(t) + \widetilde{P}(t)\widetilde{A} - \widetilde{P}(t)\widetilde{B}R^{-1}\widetilde{B}^T\widetilde{P}(t) + \widetilde{Q}^\infty = -\dot{\widetilde{P}}(t) & t \in [s, T] \\ \widetilde{P}(T) = \widetilde{Q}_f^\infty. \end{cases} \tag{57}$$

Furthermore, since each $X^N(t)^{-1}$ is continuous and well defined for each $N \in \mathbb{N}$ and that $X^\infty(t)^{-1}$ is uniformly continuous on $[s, T]$, then

$$\widetilde{P}^N(t) := Y^N(t)X^N(t)^{-1} \longrightarrow Y^\infty(t)X^\infty(t)^{-1} =: \widetilde{P}^\infty(t)$$

uniformly on $t \in [s, T]$. Finally, consider the matrix $\bar{P}(t) \in \mathbb{M}_{nM \times nM}$,

$$\bar{P}(t) = \begin{pmatrix} P(t) & 0 & \cdots & 0 \\ 0 & 0 & \cdots & 0 \\ \vdots & \vdots & \ddots & \vdots \\ 0 & 0 & \cdots & 0 \end{pmatrix},$$

where $P(t) \in \mathbb{M}_{n \times n}$ is the unique solution of the Riccati equation (54). A direct verification shows that $\bar{P}(t) \in \mathbb{M}_{nM \times nM}$ also satisfies the Riccati equation (57). However, the problem (57) admits a unique solution, implying that $\widetilde{P}^\infty(t) \equiv \bar{P}(t)$ for all $t \in [s, T]$. This concludes the proof.                                                              □

**Remark 6.4** *(Feedback optimal controls)* We would like to stress that the convergence of the Riccati matrix solution $\tilde{P}^N$ to $P$ provided by Theorem 6.3 has the following implication: For each $(s, x_0) \in [0, T] \times \mathbb{R}^n$, we define

$$\begin{aligned} \bar{u}^N(s, x_0) &:= -R^{-1}\widetilde{B}^T P^N(s)(x_0, \ldots, x_0) \text{ and} \\ \bar{u}(s, x_0) &:= -R^{-1}B^T P(s)x_0, \end{aligned} \tag{58}$$

and then, $\bar{u}^N$ tend to $\bar{u}$ for $N$ going to $+\infty$. Namely that the optimal control of problem (44), which satisfies the formula (48), evaluated at $X_0 = (x_0, \ldots, x_0)$ pointwisely converges to the optimal control of problem (42). Furthermore, for each $K \subset \mathbb{R}^n$ compact, one has that

$$\bar{u}^N(s, x_0) \to \bar{u}(s, x_0) \quad \text{uniformly on } [0, T] \times K. \tag{59}$$

It is important to point out that whereas Theorem 5.3 proves the convergence for the class of optimal open-loop controls, Theorem 6.3 deals with the convergence of optimal controls in feedback form.

It remains an open question whether Theorem 6.3 can be proved for a generic sequence of probability measures $\{\pi^N\}_{N \in \mathbb{N}}$ converging weakly-* to a generic probability measure $\pi$. Such an issue is delicate and will be studied in a forthcoming paper.

## 7 A numerical example

The aim of this section is to verify that the results summarized in Corollary 6.1 hold in a concrete example.

The model is the one presented in Sect. 6. The true dynamics is a controlled harmonic oscillator described by the matrices

$$\hat{A} := \begin{pmatrix} 0 & 1 \\ -1 & 0 \end{pmatrix} \quad \text{and} \quad B := \begin{pmatrix} 0 \\ 1 \end{pmatrix} ;$$

the coefficients that are used to define the cost functional are

$$Q := \begin{pmatrix} 1 & 0 \\ 0 & 1 \end{pmatrix} \quad \text{and} \quad R := 0.1 ,$$

and the final time is $T = 5$. Let us write $\pi^N$ as in (41) for $M = 9$, and the 9 matrices are defined in a neighborhood of matrix $\hat{A}$, i.e.,

$$A_1 = \hat{A} ,$$
$$A_{2*j+i} := \hat{A} + (-1)^i \, 0.5 \, e_j \qquad i = 0, 1, \ j = 1, 2, 3, 4$$

where $\{e_j\}_{j=1,\ldots,4}$ are the matrices of the canonical basis of $\mathbb{R}^{2 \times 2}$. The probabilities $\alpha_i^N$ are defined according to the following rule:

$$\alpha_1^N = 1 - \frac{1}{2^N} , \quad \alpha_i^N = \frac{1}{8} \frac{1}{2^N} \text{ for } i = 2, \ldots, 9 .$$

**Table 1** Errors for value functions and optimal controls related to $\pi^N$ with $N = 0, \ldots, 9$ with respect to the true value function and the true optimal control of Problem A

| $N$ | $\alpha_1$ | $\|\widetilde{V}^N - V\|_{\infty, K}$ | $Order$ | $\|\bar{u}^N - \bar{u}\|_\infty$ | $Order$ |
|---|---|---|---|---|---|
| 0 | 0 | 6.08e−0 | - | 5.25e−1 | – |
| 1 | 0.5 | 3.21e−0 | 0.92 | 3.21e−1 | 0.71 |
| 2 | 0.75 | 1.66e−0 | 0.95 | 1.82e−1 | 0.82 |
| 3 | 0.875 | 8.49e−1 | 0.97 | 9.78e−2 | 0.90 |
| 4 | 0.9375 | 4.29e−1 | 0.98 | 5.09e−2 | 0.94 |
| 5 | 0.9687 | 2.16e−1 | 0.99 | 2.59e−2 | 0.97 |
| 6 | 0.9844 | 1.08e−1 | 1.00 | 1.31e−2 | 0.99 |
| 7 | 0.9922 | 5.42e−2 | 1.00 | 6.59e−3 | 0.99 |
| 8 | 0.9961 | 2.71e−2 | 1.00 | 3.30e−3 | 1.00 |
| 9 | 0.9980 | 1.36e−3 | 1.00 | 1.65e−3 | 1.00 |

Initial point for the optimal control is $x_0 = (1, 0)$

Note that the Wasserstein distance with respect to the Euclidean norm on $\mathbb{R}^{2\times 2}$ between $\pi^N$ and $\delta_{\hat{A}}$ can be computed exactly:

$$W_1(\pi^N, \delta_{\hat{A}}) = \frac{1}{2^{N+1}} \ .$$

Both Problem A, that is the LQR problem with the matrix $\hat{A}$ known, and Problem B, that is the LQR problem with the matrix $\hat{A}$ unknown can be solved by finding the solution of a Riccati equation (see §IV.5 on Fleming–Rishel monograph [9]). We solved the equation numerically for $N = 0, \ldots, 9$. For each $N$, we computed the sup norm of the difference $\widetilde{V}^N - V$ and of the difference $\bar{u}^N - \bar{u}$, where $\bar{u}^N$ and $\bar{u}$ are, respectively, the optimal controls for the two problems starting from $x_0 = (1, 0)$. The results are summarized in Table 1. All the computations have been done using MATLAB on a MacBook Air 13" 2017 with Intel Core i5 Processor (2x 1.8 GHz).

Notice that when we increase $N$ by one, we halve the distance $W_1(\pi^N, \delta_{\hat{A}})$ and Table 1 tells us that also the error $\|\widetilde{V}^N - V\|_{\infty, K}$, with $K := [-2, 2]^2$, is halved; this is consistent with the estimate given by Corollary 6.1. At the same time, we remark that the error $\|\bar{u}^N - \bar{u}\|_\infty$ is halved as well, even if we did not have any estimate on the convergence rate of the optimal controls. We can say that in this example, both the errors are going to 0 with order 1.

The optimal trajectory of Problem A starting from $x_0 = (1, 0)$ is represented in Fig. 1. For Problem B, the optimal trajectory is actually made of 9 trajectories, the costs of which are weighted averaged in order to compute the cost functional $\widetilde{J}$. Two examples of optimal *(multi-)trajectory*, respectively, for $\pi^0$ and $\pi^2$, are represented in Figs. 2 and 3; note that the trajectory related to the true dynamics is $x^1(t)$, which is the darkest one. Finally, in Fig. 4 the optimal controls for Problem B with $N = 0, \ldots, 4$ are compared with the true optimal control.

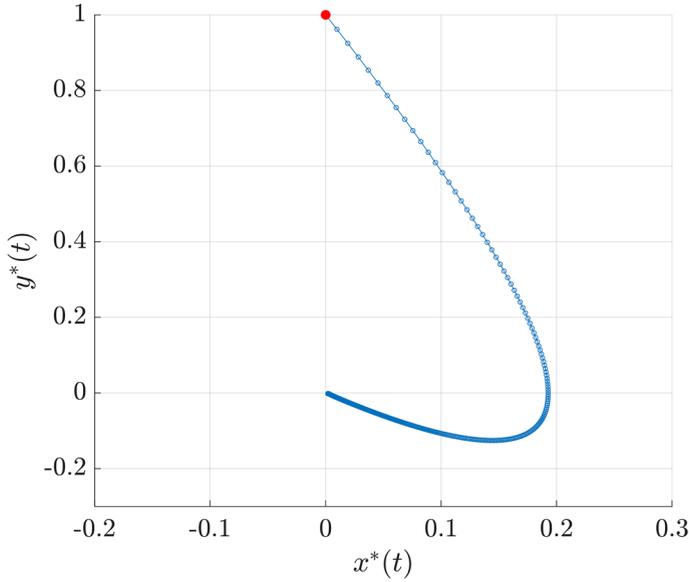

**Fig. 1** Optimal trajectory of Problem A, computed solving a 2-dimensional Riccati differential equation associated with matrix $\hat{A}$. The initial point is indicated with a red dot

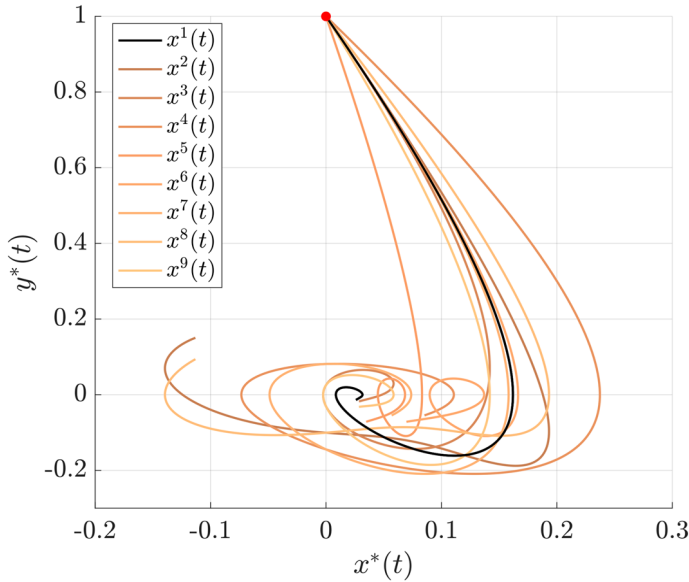

**Fig. 2** Optimal (multi-)trajectory for Problem B with $\pi^0$. To compute the optimal solution, we solved a 18-dimensional Riccati differential equation. The initial point is indicated with a red dot

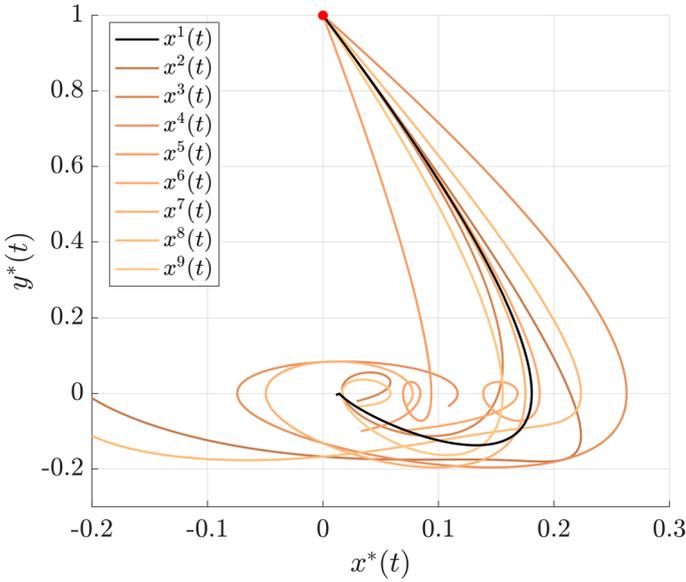

**Fig. 3** Optimal (multi-)trajectory for Problem B with $\pi^2$

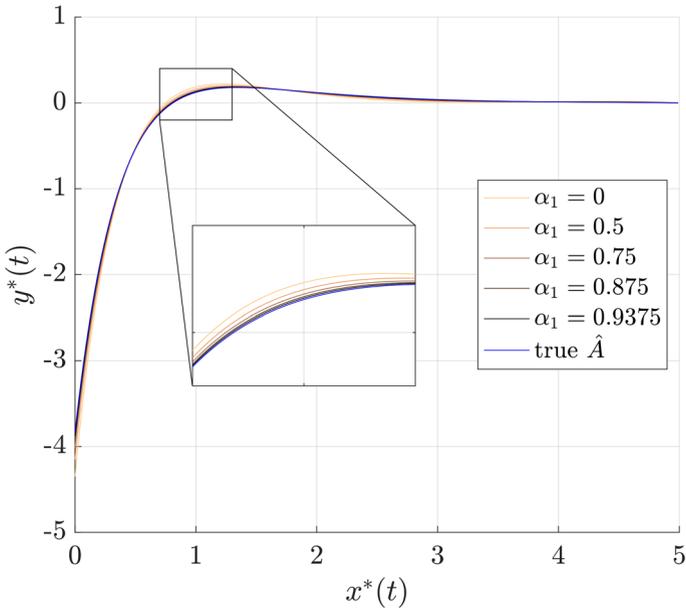

**Fig. 4** Comparison of the optimal controls of Problem B relative to different probability distributions $\pi^N$, with $N = 0, \ldots, 4$, and the true optimal control of Problem A (in blue). In the legend, we reported, for each $N$, the probability $\alpha_1$ that the true matrix $A_1 \equiv \hat{A}$ has under the distribution $\pi^N$. When $\alpha_1 \to 1$, the optimal control of Problem B converges to the true one

## 8 Conclusions

In this paper, we proved some convergence properties for the optimal policies of LQ optimal control problems with uncertainties (our Problem B), assuming that the current belief on the dynamics is represented by a generic probability distribution $\pi$ on the space of matrices. Under standard hypotheses on the system dynamics and the cost functional, we proved that the open-loop, optimal control $\bar{u}_\pi$ of Problem B converges to the open-loop, optimal control $\bar{u}$ of the actual system as soon as the distribution $\pi$ is sufficiently close (w.r.t. the Wasserstein distance (9)) to a Dirac's delta $\delta_{\hat{A}}$ evaluated at the actual system matrix $\hat{A}$. We also showed that when the probability distribution $\pi$ is actually a discrete measure, then also the closed-loop optimal control of Problem B converges to the closed-loop optimal control of the actual system when the distribution $\pi$ is sufficiently close to the Dirac's delta $\delta_{\hat{A}}$. The latter result was also validated by a numerical example presented in Sect. 7.

It is worth stressing that the proposed approach has strong connections with several Bayesian-like RL algorithms (such as PILCO), providing a theoretical framework to obtain stability and convergence guarantees for such algorithms.

As a future direction, we would like to extend this approach to a nonlinear, control affine optimal control problem with convex functional, getting closer to the problem formulation studied in [7]. A first attempt in this direction has been recently proposed in [22]. Furthermore, we are interested in constructing new efficient RL algorithms using well-established tools from control theory. In this context, we have already developed a new method for solving LQR problems with unknown dynamics [20].

**Funding** Open access funding provided by Universitá degli Studi di Roma La Sapienza within the CRUI-CARE Agreement.

**Publisher's Note**  Springer Nature remains neutral with regard to jurisdictional claims in published maps and institutional affiliations.